\documentclass{ws-sd}%
\usepackage[english]{babel}
\usepackage[utf8]{inputenc}
\usepackage[T1]{fontenc}
\usepackage{microtype}
\usepackage{enumerate}
\usepackage{mathrsfs}
\usepackage{mathtools}
\usepackage{graphicx}
\usepackage{float}
\usepackage{dsfont}
\usepackage{todonotes}
\usepackage{listings}
\usepackage{booktabs}
\usepackage[linesnumbered]{algorithm2e}
\usepackage{mathtools}
\usepackage{caption}
\usepackage{subfig}
\usepackage{tikz}
\usepackage{csquotes}
\usepackage{xcolor}
\usepackage{bbm}

\providecommand{\U}[1]{\protect\rule{.1in}{.1in}}
\newtheorem{claim}{Claim}

\newcommand\abs[1]{\left| #1 \right|}

\begin{document}

\title{ON THE MACROSCOPIC LIMIT OF BROWNIAN PARTICLES WITH LOCAL INTERACTION}
\author{FRANCO FLANDOLI}
\address{Email: franco.flandoli@sns.it.\\ Scuola Normale Superiore, Piazza dei Cavalieri, 7, 56126 Pisa, Italy.}
\author{MARTA LEOCATA}
\address{Email: leocata@math.univ-lyon1.fr. \\Institut Camille Jordan, Université Claude Bernard Lyon 1, 43 boulevard du 11 novembre 1918 F-69622 Villeurbanne Cedex, France}
\author{CRISTIANO RICCI}
\address{Email: cristiano.ricci@sns.it. \\Scuola Normale Superiore, Piazza dei Cavalieri, 7, 56126 Pisa, Italy.}
\maketitle

\begin{abstract}
An interacting particle system made of diffusion processes with local
interaction is considered and the macroscopic limit to a nonlinear PDE is
investigated. Few rigorous results exists on this problem and in particular
the explicit form of the nonlinearity is not known. The paper reviews this
subject, some of the main ideas to get the limit nonlinear PDE and provides
both heuristic and numerical informations on the precise form of the
nonlinearity which are new with respect to the literature and coherent with
the few known informations. 

\end{abstract}



\section{Introduction}

We are concerned with an interacting particle system governed by diffusions
processes as follows:%
\begin{equation}
dX_{t}^{i,N}=-\sum_{j\neq i}\nabla V\left(  X_{t}^{i,N}-X_{t}^{j,N}\right)
dt+\sigma dB_{t}^{i} \label{model in large box}%
\end{equation}
where $i=1,...,K_{N}$, $X_{t}^{i,N}\in\mathbb{R}^{d}$, $V:\mathbb{R}%
^{d}\rightarrow\mathbb{R}$ has the form%
\begin{equation}
V\left(  x\right)  =U\left(  \left\vert x\right\vert \right)
\label{potential invariant}%
\end{equation}
where $U:\left(  0,\infty\right)  \rightarrow\mathbb{R}$ is twice
differentiable, either with compact support or a suitable decay at infinity,
$B_{t}^{i}$ are independent Brownian motions in $\mathbb{R}^{d}$ and
$\sigma>0$.

The number $N$ is the scaling parameter and the number of particles $K_{N}$
will be assumed of order $N^{d}$. For mathematical simplicity we assume the
particles live in a large torus%
\[
\mathbb{T}_{N}^{d}:=\mathbb{R}^{d}/N\mathbb{Z}^{d}%
\]
(the set $\left[  0,N\right]  ^{d}$ with periodic identifications). When we
assume
\[
K_{N}=\left\lfloor \rho N^{d}\right\rfloor
\]
for some $\rho>0$, we are saying that the density of particles in
$\mathbb{T}_{N}^{d}$ is $\rho$.

We assume that the initial conditions $X_{0}^{i,N}$, $i=1,...,N^{d}$, are
random, independent, with a distributions such that the typical distance
between neighbor particles is of order one, or $\frac{1}{\rho}$, but not
concentrated with infinitesimal-in-$N$ interparticle distance. Expecting the
same holds for $t>0$, each particle $X_{t}^{i,N}$ typically interacts only
with a finite number of other particles in the case when $K$ is compact
support (by finite number we mean finite in the limit when $N$ goes to
infinity). Or, when $U$ has infinite support but decays suitably at infinity,
although the number of particles seen by $X_{t}^{i,N}$ is infinite, only a
finite number has a relevant influence on $X_{t}^{i,N}$. This is \textit{not}
a mean field regime; we call it local interaction regime. We shall also
comment on intermediate situations between the two.

We want to investigate the macroscopic limit of this system, namely the weak
limit of the empirical measure%
\[
\mu_{t}^{N}=\frac{1}{N^{d}}\sum_{i=1}^{K_{N}}\delta_{\frac{1}{N}X_{t\cdot
N^{2}}^{i,N}}%
\]
corresponding to a parabolic zoom in space and time, natural because the
transformation $B_{t}^{i}\rightarrow\frac{1}{N}B_{t\cdot N^{2}}^{i}$ leaves
the law of Brownian motion invariant. Notice that $\mu_{t}^{N}$ is not a
probability measure, unless $K_{N}=N^{d}$.

In the local interaction case considered here this limit is still poorly
understood. In the works \cite{varadhan1991scaling}, \cite{uchiyama2000pressure}, similarly to what happens for hydrodynamic limits of discrete systems \cite{kipnis1998scaling}, it is proved that the weak limit of the empirical measure $\mu^N_t$ is a weak solution of the following nonlinear Partial Differential Equation (PDE):
\[\partial_t\rho=\frac{1}{2}\Delta P_V(\rho).\]
but, apart from a number of restrictions on V imposed in these works, the main gap with respect to the discrete case is the lack of quantitative information on $P_V(\rho)$.
In this exploratory work we review some facts known in the
literature and present conjectures based on heuristic arguments and numerical
simulations. We distinguish between the case when the interaction is purely
repulsive, namely $U^{\prime}\left(  r\right)  <0$ for $r>0$ (possibly only up
to some $r_{1}>0$ beyond which $U^{\prime}\left(  r\right)  =0$) and the case
when the interaction includes an attractive component, namely $U^{\prime
}\left(  r\right)  <0$ for $0<r<r_{0}$ and $U^{\prime}\left(  r\right)  >0$
for $r>r_{0}$ (again possibly only up to some $r_{1}>r_{0}$ beyond which
$U^{\prime}\left(  r\right)  =0$). The repulsive case is better
understood;\ the case with also local attraction is very difficult, with
several obscure aspects.

Our motivation for studying this problem has been the desire to model
\textit{adhesion} between cells. A possible way of modeling such phenomenon is given by hard-core interacting particles. In this way particles are thought to be hard spheres which cannot compenetrate at all. Between the first results on the continuum limit for this type of particle system, we mention for one dimension  \cite{rost1984diffusion}  and for higher dimension \cite{bruna2012excluded}. A very recent result in this direction is \cite{gavish2019large}. The kind of interaction in which we are interested is the different from the one above mentioned.
Assume a family of living cells is modeled
simply by a position $X_{t}^{i,N}$ and a local interaction. Repulsion is
motivated by a volume constraint: a cell is not a point, it has a finite size,
and called $r_{0}$ its diameter, repulsion acts when the centers of the cells
are at a distance smaller than $r_{0}$. But when they are at a distance
slightly larger than $r_{0}$, cells do not simply separate: they are linked by
macromolecules that produce adhesion between the cell membranes. We may steer
the distance between the centers up so some value without splitting the cells,
which are then subject to an attractive force. After some distance, the cells
separate and do not feel each other anymore, corresponding to a compact
support function $U$. 
Most of the literature describes cell adhesion by
non-local attracting forces which are not realistic; they are a simplification
(since they lead to mean field theories) and may give relatively good
quantitative results when the mean field kernel has very short range, see \cite{armstrong2006continuum},\cite{buttenschoen2018space} \cite{flandoli_leocata_2019};\ but
conceptually these models are wrong, since each cell interacts with infinitely
many others and in a weak uniform way. The literature in biomathematics on
cell adhesion seems to ignore the possibility, offered by the works of
Varadhan \cite{varadhan1991scaling} and Uchiyama \cite{uchiyama2000pressure}, of studying the macroscopic limit of system
like (\ref{model in large box}) having true local interaction. Our motivation
for writing this work is to popularize this bibliographical link and propose
additional quantitative conjectures beyond those made in \cite{varadhan1991scaling}, \cite{uchiyama2000pressure}.

\section{Macroscopic view}

We now zoom and observe the previous particles as they were very close points
in the unitary torus $\mathbb{T}^{d}=\mathbb{R}^{d}/\mathbb{Z}^{d}$ and we
accelerate time according to the invariance of Brownian motion;\ we introduce
the notations:%
\begin{align*}
Y_{t}^{i,N}  &  :=\frac{1}{N}X_{t\cdot N^{2}}^{i,N}\in\mathbb{T}^{d}\\
W_{t}^{i}  &  :=\frac{1}{N}B_{t\cdot N^{2}}^{i}%
\end{align*}
recalling that $W_{t}^{i}$ are independent Brownian motions. We have%
\[
dY_{t}^{i,N}=-N\sum_{j\neq i}\nabla V\left(  X_{t\cdot N^{2}}^{i,N}-X_{t\cdot
N^{2}}^{j,N}\right)  dt+\sigma dW_{t}^{i}.
\]
Set%
\[
V_{N}\left(  x\right)  :=N^{d}V\left(  Nx\right)  =N^{d}U\left(  \left\vert
Nx\right\vert \right)  .
\]
Then
\[
\frac{1}{N^{d}}\nabla V_{N}\left(  x\right)  =N\nabla V\left(  Nx\right)
\]
and thus we may write the previous rescaled equation in the elegant form%
\begin{equation}
dY_{t}^{i,N}=-\frac{1}{N^{d}}\sum_{j\neq i}\nabla V_{N}\left(  Y_{t}%
^{i,N}-Y_{t}^{j,N}\right)  dt+\sigma dW_{t}^{i} \label{mean field like}%
\end{equation}
which sounds like a mean field equation, due to the factor $\frac{1}{N^{d}}$,
but it is not because the potential is rescaled (and it has infinitesimal
range of interaction, in $\mathbb{T}^{d}$).

By It\^{o} formula, if $\phi:\mathbb{T}^{d}\rightarrow\mathbb{R}$ is a smooth
compact support test function, then
\begin{align*}
d\phi\left(  Y_{t}^{i,N}\right)   &  =-\left(  \nabla\phi\right)  \left(
Y_{t}^{i,N}\right)  \frac{1}{N^{d}}\sum_{j\neq i}\nabla V_{N}\left(
Y_{t}^{i,N}-Y_{t}^{j,N}\right)  dt\\
&  +\left(  \nabla\phi\right)  \left(  Y_{t}^{i,N}\right)  \sigma dW_{t}%
^{i}+\frac{\sigma^{2}}{2}\Delta\phi\left(  Y_{t}^{i,N}\right)  dt.
\end{align*}
Let us use a notational trick: the function $U$ is not defined for $r=0$,
hence $V_{N}\left(  0\right)  $ is not defined; we set it equal to
zero\footnote{By this definition, $\nabla V_{N}\left(  Y_{t}^{i,N}-Y_{t}%
^{j,N}\right)  =0$ when $j=i$. However, more subtle is the problem that we
could have $Y_{t}^{i,N}-Y_{t}^{j,N}=0$ also sometimes for $j\neq i$.
Fortunately, since we always assume to have a repulsive component in the
interaction,  one can prove this never happens
.}. Hence we may remove
the restriction $j\neq i$ in the sum and write%
\begin{align*}
d\phi\left(  Y_{t}^{i,N}\right)   &  =-\left(  \nabla\phi\right)  \left(
Y_{t}^{i,N}\right)  \frac{1}{N^{d}}\sum_{j=1}^{K_{N}}\nabla V_{N}\left(
Y_{t}^{i,N}-Y_{t}^{j,N}\right)  dt+dR_{t}^{N,\phi}\\
&  =-\left(  \nabla\phi\right)  \left(  Y_{t}^{i,N}\right)  \int%
_{\mathbb{T}^{d}}\nabla V_{N}\left(  Y_{t}^{i,N}-y\right)  \mu_{t}^{N}\left(
dy\right)  dt+dR_{t}^{N,\phi}%
\end{align*}
where $dR_{t}^{N,\phi}=\left(  \nabla\phi\right)  \left(  Y_{t}^{i,N}\right)
\sigma dW_{t}^{i}+\frac{\sigma^{2}}{2}\Delta\phi\left(  Y_{t}^{i,N}\right)
dt$. Therefore%
\[
d\left\langle \phi,\mu_{t}^{N}\right\rangle =-\left\langle \nabla\phi
\int_{\mathbb{T}^{d}}\nabla V_{N}\left(  \cdot-y\right)  \mu_{t}^{N}\left(
dy\right)  ,\mu_{t}^{N}\right\rangle dt+\frac{\sigma^{2}}{2}\left\langle
\Delta\phi,\mu_{t}^{N}\right\rangle dt+dM_{t}^{N,\phi}%
\]
where the martingale $M_{t}^{N,\phi}$ is given by
\[
M_{t}^{N,\phi}=\int_{0}^{t}\frac{1}{N^{d}}\sum_{j=1}^{K_{N}}\left(  \nabla
\phi\right)  \left(  Y_{s}^{i,N}\right)  \sigma dW_{s}^{\left(  i\right)  }.
\]
A classical simple computation with the isometry formula of It\^{o} calculus
proves that $M_{t}^{N,\phi}$ converges to zero in mean square; and also
uniformly in time, using Doob's inequality. Assuming one can prove that
$\mu_{t}^{N}$ converges weakly to a measure $\mu_{t}$, uniformly in time
(maybe up to subsequences) and that $\mu_{t}$ has density $\rho_{t}$ with
respect to Lebesgue measure (this is not necessary immediately, but will play
a role later on), under the assumption that the same holds for the initial
condition, we may pass to the limit in the terms $\left\langle \phi,\mu
_{t}^{N}\right\rangle $, $\left\langle \phi,\mu_{0}^{N}\right\rangle $,
$\int_{0}^{t}\left\langle \Delta\phi,\mu_{s}^{N}\right\rangle ds$,
$M_{t}^{N,\phi}$;\ hence also the remaining term has a limit and we get%
\begin{equation}
\left\langle \phi,\rho_{t}\right\rangle =\left\langle \phi,\rho_{0}%
\right\rangle -\lim_{N\rightarrow\infty}\int_{0}^{t}\left\langle \nabla
\phi\int_{\mathbb{T}^{d}}\nabla V_{N}\left(  \cdot-y\right)  \mu_{s}%
^{N}\left(  dy\right)  ,\mu_{s}^{N}\right\rangle ds+\frac{\sigma^{2}}{2}%
\int_{0}^{t}\left\langle \Delta\phi,\rho_{s}\right\rangle
ds.\label{PDE implicit}%
\end{equation}
The main problem is to identify the limit left implicit above. For the purpose
of the overview, we first identify the limit in the classical mean field case,
Section \ref{Sect mean field}; then we identify it in the case of repulsive
integrable potential by means of simplified arguments not properly of local
type, Sections \ref{sect two steps}, \ref{sect intermediate}, and finally we
discuss it in the main case motivating this paper, namely the case of local
interaction, Section \ref{Sect local}.

The tightness of the family of laws of $\mu^{N}$ required to implement
rigorously the previous arguments is not trivial and it is discussed in \cite{varadhan1991scaling},
\cite{uchiyama2000pressure}, under different conditions.

\subsection{Mean field interaction\label{Sect mean field}}

In broad terms, the mean field case is when we start directly in the unitary
torus $\mathbb{T}^{d}$ with equations (\ref{mean field like}) but with $V_{N}$
independent of $N$:%
\begin{equation}
V_{N}\left(  x\right)  =V\left(  x\right)  \text{.}%
\label{potential invariant 2}%
\end{equation}
Going back to formulation (\ref{model in large box}) in the large box
$\mathbb{T}_{N}^{d}$, the potential there should depend on $N$. Thus the mean
field case is not a particular case of the problem studied in this paper,
since we started from (\ref{model in large box}) with a given potential. The
only particular case satisfying (\ref{potential invariant}),
(\ref{potential invariant 2}) and $V_{N}\left(  x\right)  :=N^{d}V\left(
Nx\right)  $ is the case
\begin{equation}
V\left(  x\right)  =\frac{1}{\left\Vert x\right\Vert ^{d}}%
.\label{mean field conditon}%
\end{equation}
As a curiosity, this is the boundary case between weak and strong repulsion
described in Section \ref{sect remarks potential}.

When (\ref{potential invariant 2}) is imposed in equation
(\ref{mean field like}), and $\nabla V$ is continuous and bounded, then
tightness of the family of laws of $\mu^{n}$ is much easier [Sznitman] and
convergence of the nonlinear term in (\ref{PDE implicit}) is almost trivial:%
\[
\int_{\mathbb{T}^{d}}\nabla V_{N}\left(  \cdot-y\right)  \mu_{s}^{N}\left(
dy\right)  \rightarrow\left\langle \nabla V\left(  \cdot-y\right)  ,\rho
_{s}\right\rangle =:\left(  \nabla V\ast\rho_{s}\right)  \left(  y\right)
\]
and%
\[
\int_{0}^{t}\left\langle \nabla\phi\int_{\mathbb{T}^{d}}\nabla V_{N}\left(
\cdot-y\right)  \mu_{s}^{N}\left(  dy\right)  ,\mu_{s}^{N}\right\rangle
ds\rightarrow\int_{0}^{t}\left\langle \nabla\phi\cdot\left(  \nabla V\ast
\rho_{s}\right)  ,\rho_{s}\right\rangle ds.
\]
Integrating (formally)\ by parts we get the mean field equation%
\[
\partial_{t}\rho=\frac{\sigma^{2}}{2}\Delta\rho+\operatorname{div}\left(
\rho\left(  \nabla V\ast\rho\right)  \right)  .
\]
In applications, taking $V$ with very small support is a practical way to get
numerical simulations very close to adhesion. But obviously the model required
long range interaction, so it is logically incorrect, although reasonable
under the rough view of a numerical simulation.

\subsection{Two-step limit under integrable repulsive
potential\label{sect two steps}}

Although not logically correct, there is a cheap way to obtain a guess about
the limit in (\ref{PDE implicit}). It is based on a limit taken in two
successive steps. We introduce two scaling parameters $N$ and $M$ and replace
the limit in (\ref{PDE implicit}) by%
\[
\lim_{M\rightarrow\infty}\lim_{N\rightarrow\infty}\int_{0}^{t}\left\langle
\nabla\phi\int_{\mathbb{T}^{d}}\nabla V_{M}\left(  \cdot-y\right)  \mu_{s}%
^{N}\left(  dy\right)  ,\mu_{s}^{N}\right\rangle ds.
\]
The first limit, in $N$, is like the mean field case (under the assumption
that $\nabla V_{M}$ is continuous and bounded) and thus we get%
\[
\lim_{M\rightarrow\infty}\int_{0}^{t}\left\langle \nabla\phi\cdot\left(
\nabla V_{M}\ast\rho_{s}^{M}\right)  ,\rho_{s}^{M}\right\rangle ds.
\]
Now assume a particular but natural version of the \textit{repulsive}
case:\ assume that $V$ is a probability density, of the form $V\left(
x\right)  =U\left(  \left\vert x\right\vert \right)  $ with $U$ decreasing on
$\left(  0,\infty\right)  $. More precisely, assume it is such after
normalization by
\[
C_{V}:=\int V\left(  x\right)  dx.
\]
Then $C_{V}^{-1}V_{M}\left(  x\right)  :=M^{d}C_{V}^{-1}V\left(  Mx\right)  $
are classical mollifiers, with the property that
\[
\int_{\mathbb{T}^{d}}C_{V}^{-1}V_{M}\left(  x-y\right)  f\left(  x\right)
dx\rightarrow f\left(  y\right)
\]
where convergence is for instance uniform on bounded sets when $f$ is
uniformly continuous (several other results are known under different
assumptions on $f$). Hence, assuming we can prove that $\rho_{s}^{M}$
converges to a limit $\rho_{s}$ in a suitable topology compatible with results
of convergence of mollifiers (maybe up to subsequences), we have%
\begin{align*}
\left(  \nabla V_{M}\ast\rho_{s}^{M}\right)  \left(  y\right)   &
=\int_{\mathbb{T}^{d}}\nabla V_{M}\left(  x-y\right)  \rho_{s}^{M}\left(
x\right)  dx=-\int_{\mathbb{T}^{d}}V_{M}\left(  x-y\right)  \nabla\rho_{s}%
^{M}\left(  x\right)  dx\\
&  \rightarrow-C_{V}\nabla\rho_{s}\left(  y\right)  \qquad\text{as
}M\rightarrow\infty
\end{align*}
and thus (up to rigorous care)%
\begin{align*}
\lim_{M\rightarrow\infty}\int_{0}^{t}\left\langle \nabla\phi\cdot\left(
\nabla V_{M}\ast\rho_{s}^{M}\right)  ,\rho_{s}^{M}\right\rangle ds &
=-C_{V}\int_{0}^{t}\left\langle \nabla\phi\cdot\nabla\rho_{s},\rho
_{s}\right\rangle ds\\
&  =-\frac{C_{V}}{2}\int_{0}^{t}\left\langle \nabla\phi,\nabla\rho_{s}%
^{2}\right\rangle ds\\
&  =\frac{C_{V}}{2}\int_{0}^{t}\left\langle \Delta\phi,\rho_{s}^{2}%
\right\rangle ds.
\end{align*}
Here we see for the first time the role of the density $\rho_{s}$ with respect
to the measure $\mu_{s}\left(  x\right)  =\rho_{s}\left(  x\right)  dx$:\ we
need to take the square $\rho_{s}^{2}$, meaningful only for densities. The
final equation, formally written after integration by parts, is%
\begin{equation}
\partial_{t}\rho=\frac{\sigma^{2}}{2}\Delta\rho+\frac{C_{V}}{2}\Delta\rho
^{2}.\label{intermediate PDE}%
\end{equation}

Notice that this result required $V$ integrable, $C_{V}<\infty$. Under this
assumption, a result of Lemma 8.5 of \cite{uchiyama2000pressure} plus (8.10) of the same paper prove
rigorously that the large $\rho$ asymptotic of the nonlinearity in
(\ref{intermediate PDE}) is precisely $\rho^{2}$. Our numerical simulations
for repulsive integrable potentials confirm (\ref{intermediate PDE}) also from
other quantitative sides (not only the degree two for large $\rho$). Equation
(\ref{intermediate PDE}), in all its quantitative aspects, is also rigorously
proved in the intermediate regime described in the next section. Therefore it
seems that the simple conjecture based on the two-step method is quite
realistic. 

\begin{figure}[t]
\centering
\begin{minipage}[b]{0.4\textwidth}
\includegraphics[width=\textwidth]{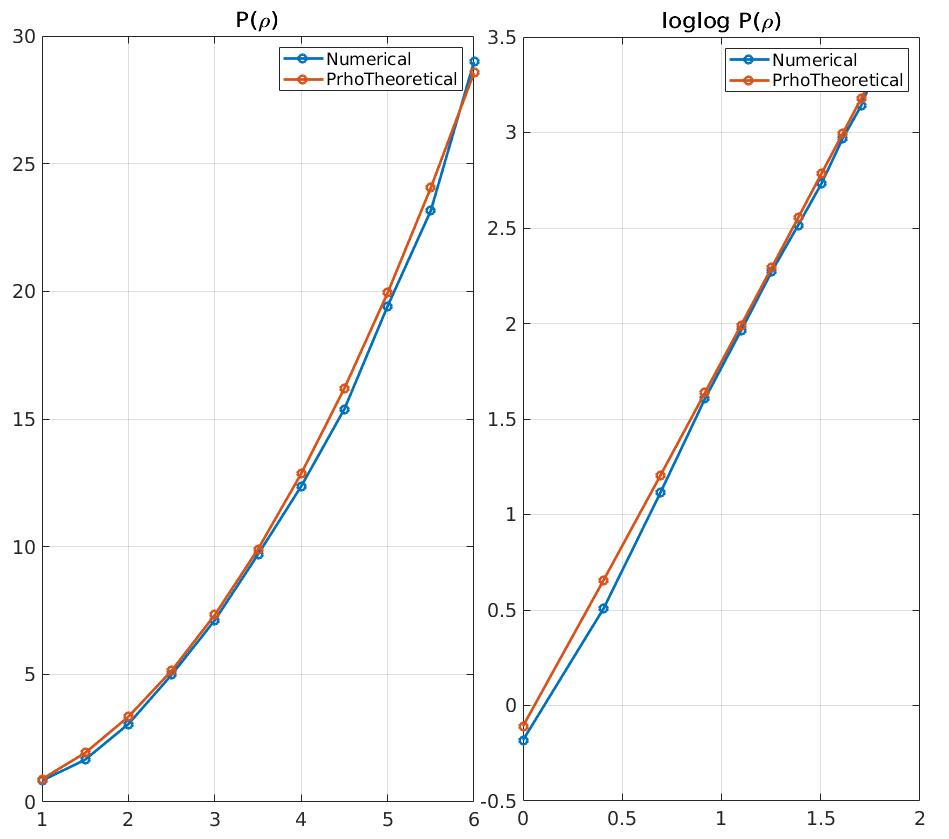}
\end{minipage}\hspace{1.5cm}
\begin{minipage}[t]{0.3\textwidth}
\vspace{-6.5cm}
\begin{tikzpicture}[domain=0:4,scale=0.5]
    \draw[->] (-0.5,0) -- (5,0) ; 
    \draw[->] (0,-0.5) -- (0,6);
	\draw[color=red,thick] plot (\x, 5*exp(-\x^2/0.9);
\end{tikzpicture} 

\begin{tikzpicture}[domain=0:4,scale=0.5]
    \draw[->] (-0.5,0) -- (5,0) ; 
    \draw[->] (0,-0.5) -- (0,6);
	\draw[color=red,thick] plot ({\x}, {5*(1-\x^3)*exp(-\x^2/0.9)});
\end{tikzpicture} 
\end{minipage}
\caption{Left: Function $P_{V}(\rho)$ in the case of potential \eqref{eq:attrgauus}. Comparison between the function $P_{V}(\rho) $ obtained by numerical simulations and the function $ \sigma^{2}\rho + C_{V}\rho^{2}$. Left: comparison between functions $P_{V}(\rho)$ in natural scale. Center: Comparison in $\log\log$ scale.  Right: plot of the potentials in dimension one $V(x)$ equal to \eqref{eq:repgauss} (top) and \eqref{eq:attrgauus} (bottom).}\label{fig:attrGauss}
\end{figure}
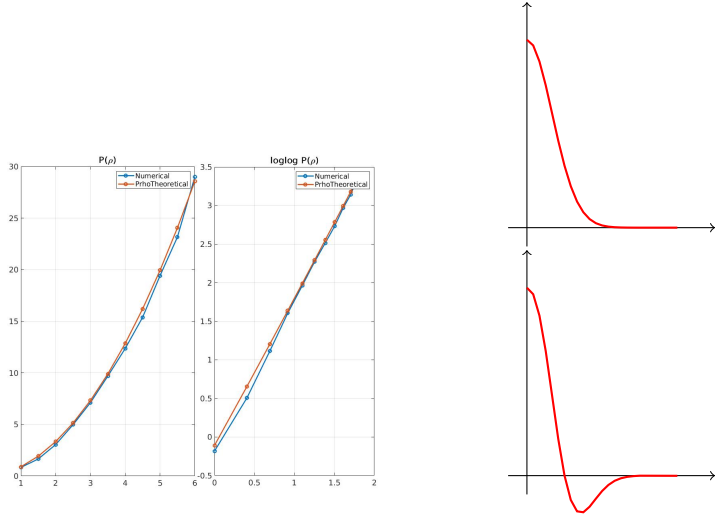

In order to validate this conjecture we present here some numerical results. Here we only briefly present the results, since the numerical analysis behind  is more carefully explained in Section \ref{sec:quant}. We consider the following two potentials, whose plot is represented in Figure \ref{fig:attrGauss} (right):
\begin{equation}\label{eq:repgauss}
V(x) = \exp\left(-\frac{\abs{x}^{2}}{2\sigma^{2}}\right)	
\end{equation}
and
\begin{equation}\label{eq:attrgauus}
V(x) = \exp\left(-\frac{\abs{x}^{2}}{2\sigma^{2}}\right)(1-x^{3}).	
\end{equation}
In Figure \ref{fig:attrGauss} (left) we superimpose the function $P_{V}(\rho) = \sigma^{2}\rho+C_{V}\rho^{2} $, for $V$ equal to \eqref{eq:attrgauus}, with the function $P_{V}(\rho)$ computed numerically. As we can see there is almost perfect superimposition of the two, confirming our conjecture also for attractive potentials with negative sign, like \eqref{eq:attrgauus}. We also remark that simulations confirm that the leading term in the function $P_{V}(\rho)$  is of order $\rho^{2}$. It is however important to notice that Figure \ref{fig:attrGauss} contains information only for $\rho > 1$. For smaller values of $\rho$ in fact there are some discrepancy between numerics and theory that are still not fully understood.
For potential \eqref{eq:repgauss} we mention that the function $P_{V}(\rho) = \sigma^{2}\rho+C_{V}\rho^{2} $
and that computed numerically coincide also for values of $\rho$ smaller than one, confirming fully equation \eqref{intermediate PDE} as limiting equation.
We also remark that a similar result holds, with the same constant $C_{V}$, in
Bose-Einstein condensation theory, see \cite[Chapter 7]{rougerie2015finetti}.

\subsection{Intermediate interactions\label{sect intermediate}}

Karl Oelschläger, in a series of papers \cite{oelschlager1985law}, \cite{oelschlager1989derivation},\cite{oelschlager1990large}, clarified rigorously the
results when the interaction is intermediate between purely local and mean
field. It is assumed that
\[
V_{N}\left(  x\right)  :=N^{\beta d}V\left(  Nx\right)
\]
for some%
\[
0<\beta<1.
\]
The case $\beta=1$ corresponds to local interaction, $\beta=0$ to mean field
interaction. Intuitively, each particle $X_{t}^{i,N}$ interacts with
infinitely many others but still with an infinitely small proportion of the total.

The potential $V$, in those works, is repulsive and integrable, as in the
two-step approach described above. The final result is always equation
(\ref{intermediate PDE}). This is the best available confirmation that
(\ref{intermediate PDE}) is the correct one in the repulsive integrable regime
and, as already said, our simulations confirm the result also in the limit
case of local interactions.

\section{Local Interaction\label{Sect local}}

When the interaction potential $V$ is not integrable or not just repulsive,
equation (\ref{intermediate PDE}) seems to be false; it seems it should be
replaced by an equation of the form%
\begin{equation}
\partial_{t}\rho=\frac{1}{2}\Delta P_{V}\left(  \rho\right)  \label{PDE}%
\end{equation}
where the nonlinear function $P_{V}\left(  \rho\right)  $ is close to linear
for small values of $\rho$ (like $\sigma^{2}\Delta\rho+C_{V}\Delta\rho^{2}$)
but growth more than quadratically as $\rho\rightarrow\infty$ and, when there
is also suitable attraction, the slope for small $\rho$ may be different from
$\sigma^{2}$. The first rigorous proof that (\ref{PDE}) is the macroscopic
limit of system (\ref{model in large box}) has been given in $d=1$ by Varadhan
\cite{varadhan1991scaling} in the case of repulsive potential, without a quantification of
$P_{V}\left(  \rho\right)  $.\ In general dimension, the basic result has been
proved by Uchiyama \cite{uchiyama2000pressure};\ it is a conditional result (conditional because it
is based on an ergodic assumption which is an open problem) but holds for
quite general potentials, which may include also an attractive part; and in
the Appendix, in particular Lemma 8.5 (see also (8.10)), \cite{uchiyama2000pressure} provides
quantitative result of the form $P_{V}\left(  \rho\right)  =O\left(
\rho^{\gamma}\right)  $ as $\rho\rightarrow\infty$ with precise prescription
of $\gamma$ depending on the singularity of $V$ at $x=0$. Thus, although being
only a one-side result (because it is of the form $O\left(  \rho^{\gamma
}\right)  $) and it is only for large $\rho$, it is a precious indications.
Not only it is the only quantitative rigorous result but also the order
$\gamma$ of the power is confirmed by our numerical simulations to be the true
one, not only a bound on one side.

Below, our goal is to give more explicit forms of $P_{V}\left(  \rho\right)
$, coherent with \cite{uchiyama2000pressure}, obtained when possible by both heuristic arguments and
numerical simulations.

\subsection{Some remark on the potential\label{sect remarks potential}}

As already remarked above, we have in mind sometimes the purely repulsive
case, since it is easier, and sometime else the attractive-repulsive case,
motivated for instance by cell adhesion. Let us describe the
attractive-repulsive case. The potential $U$ is made of two components, a
repulsive one corresponding to volume constraint $V_{vol}$ and an attractive
one corresponding to adhesion $V_{adh}$
\[
U(r)=U_{vol}(r)+U_{adh}(r).
\]
where $U,U_{adh},U_{vol}:(0,+\infty)\rightarrow\mathbb{R}$ are piecewise
$C^{1}$ functions with the properties%
\[
U_{vol}^{\prime}%
\begin{cases}
<0\quad\text{if }r\in(0,r_{1})\\
=0\quad\text{if }r>r_{1}%
\end{cases}
\]
and
\[
U_{adh}^{\prime}%
\begin{cases}
<0\quad\text{if }r\in(0,r_{0})\\
>0\quad\text{if }r\in(r_{0},r_{1})\\
=0\quad\text{if }r>R_{1}.
\end{cases}
\]
This is the general form in the compact support case, but we shall also
consider the case when both $U_{adh}$ and $U_{vol}$ have full support but
suitable decay at infinity, like the case of Lennard-Jones potential.

As mentioned several times above, the most relevant results on macroscopic
limit for local interaction are given by \cite{varadhan1991scaling} and
\cite{uchiyama2000pressure}. In both paper the result of convergence is stated
under some condition on the potential, see for instance page 1160 of
\cite{uchiyama2000pressure}, which include both the case of a purely repulsive
potential $U(r)=U_{vol}(r)$ and the case of attractive-repulsive one
$U(r)=U_{vol}(r)+U_{adh}(r)$. For certain results a main distinction is made
by the following condition:\ we call \textit{weak repulsion} the case when
$\int_{0}^{1}U(r)r^{d-1}dr<\infty$, \textit{strong repulsion} the case when
$\int_{0}^{1}U(r)r^{d-1}dr=\infty$. 

\subsection{Manipulation of the nonlinear term}

Apart from the tightness problem that we address to the literature, the main
problem left open above is to understand the limit of nonlinear term in
equation (\ref{PDE implicit})%
\[
\left\langle \mu_{t}^{N},\nabla\phi\cdot\nabla\left(  V_{N}\ast\mu_{t}%
^{N}\right)  \right\rangle =\left\langle \nabla\phi\int_{\mathbb{T}^{d}}\nabla
V_{N}\left(  \cdot-y\right)  \mu_{t}^{N}\left(  dy\right)  ,\mu_{t}%
^{N}\right\rangle .
\]
Let us assume compact support potential, being $r_{1}$ the size of the support
as in the examples above. Let us rewrite the nonlinear term explicitly as
(assuming also $K_{N}=\left[  N^{d}\rho\right]  $)%
\[
=\frac{1}{N^{d}}\sum_{i=1}^{\left[  N^{d}\rho\right]  }\nabla\phi\left(
Y_{t}^{i,N}\right)  \cdot N\sum_{j:\left\vert Y_{t}^{i,N}-Y_{t}^{j,N}%
\right\vert \leq\frac{r_{1}}{N}}U^{\prime}\left(  N\left\vert Y_{t}%
^{i,N}-Y_{t}^{j,N}\right\vert \right)  \frac{Y_{t}^{i,N}-Y_{t}^{j,N}%
}{\left\vert Y_{t}^{i,N}-Y_{t}^{j,N}\right\vert }%
\]
or also as%
\[
=\frac{1}{N^{2d}}\sum_{i,j=1}^{\left[  N^{d}\rho\right]  }\nabla\phi\left(
Y_{t}^{i,N}\right)  \cdot\nabla V_{N}\left(  Y_{t}^{i,N}-Y_{t}^{j,N}\right)  .
\]
We shall use both expressions.

There is a cancellation we need to implement;\ systems with this cancellation
are called gradient systems in other contexts \cite{kipnis1998scaling};\ and the same
cancellation is used in Schoquet symmetrization approach to measure valued
solutions to 2D Euler equations \cite{schochet1995weak}. Under our assumptions we have%
\[
\nabla V_{N}\left(  -x\right)  =-\nabla V_{N}\left(  x\right)  .
\]
Observe two particles $X_{t}^{i}$ and $X_{t}^{j}$: the ``force'' impressed by
$X_{t}^{j}$ on $X_{t}^{i}$ is
\[
-\frac{1}{N^{d}}\nabla V_{N}\left(  Y_{t}^{i,N}-Y_{t}^{j,N}\right)
\qquad\text{(force of }{\small X}_{t}^{j}\text{ on }{\small X}_{t}^{i}\text{)}%
\]
while the force impressed by $X_{t}^{i}$ on $X_{t}^{j}$ is
\[
-\frac{1}{N^{d}}\nabla V_{N}\left(  Y_{t}^{i,N}-Y_{t}^{j,N}\right)  =\frac
{1}{N^{d}}\nabla V_{N}\left(  Y_{t}^{i,N}-Y_{t}^{j,N}\right)  \qquad
\text{(force of }Y_{t}^{j,N}\text{ on }Y_{t}^{i,N}\text{).}%
\]
In the equations of motion these two forces never appear together, but in the
formula for $\left\langle S_{t}^{N},\nabla\phi\cdot\nabla\left(  V_{N}\ast
S_{t}^{N}\right)  \right\rangle $ the analogous terms appear together,
precisely as:%
\begin{align*}
&  \nabla\phi\left(  Y_{t}^{i,N}\right)  \cdot\nabla V_{N}\left(  Y_{t}%
^{i,N}-Y_{t}^{j,N}\right)  +\nabla\phi\left(  Y_{t}^{j,N}\right)  \cdot\nabla
V_{N}\left(  Y_{t}^{j,N}-Y_{t}^{i,N}\right) \\
&  =\left(  \nabla\phi\left(  Y_{t}^{i,N}\right)  -\nabla\phi\left(
Y_{t}^{j,N}\right)  \right)  \cdot\nabla V_{N}\left(  Y_{t}^{i,N}-Y_{t}%
^{j,N}\right)  .
\end{align*}
Since $\phi$ is smooth, and only nearby particles interact, this is almost a
cancellation. By Taylor formula it is approximatively equal to
\begin{align*}
&  =\sum_{\alpha=1}^{d}\partial_{\alpha}\left(  \phi\left(  Y_{t}%
^{i,N}\right)  -\phi\left(  Y_{t}^{j,N}\right)  \right)  \partial_{\alpha
}V_{N}\left(  Y_{t}^{i,N}-Y_{t}^{j,N}\right) \\
&  \sim\sum_{\alpha,\beta=1}^{d}\partial_{\alpha}\partial_{\beta}\phi\left(
Y_{t}^{j,N}\right)  \partial_{\alpha}V_{N}\left(  Y_{t}^{i,N}-Y_{t}%
^{j,N}\right)  \left(  Y_{t}^{i,N}-Y_{t}^{j,N}\right)  _{\beta}.
\end{align*}
The approximation is reasonable because $\left\vert Y_{t}^{i,N}-Y_{t}%
^{j,N}\right\vert \leq\frac{r_{1}}{N}$.

Let us put this approximation in the full nonlinear expression:%
\begin{align*}
\left\langle S_{t}^{N},\nabla\phi\cdot\nabla\left(  V_{N}\ast S_{t}%
^{N}\right)  \right\rangle  &  =\frac{1}{N^{2d}}\sum_{i<j}\sum_{\alpha
,\beta=1}^{d}\partial_{\alpha}\partial_{\beta}\phi\left(  Y_{t}^{j,N}\right)
\partial_{\alpha}V_{N}\left(  Y_{t}^{i,N}-Y_{t}^{j,N}\right)  \left(
Y_{t}^{i,N}-Y_{t}^{j,N}\right)  _{\beta}\\
&  =\frac{1}{2N^{2d}}\sum_{i,j=1}^{\left[  N^{d}\rho\right]  }\sum
_{\alpha,\beta=1}^{d}\partial_{\alpha}\partial_{\beta}\phi\left(  Y_{t}%
^{j,N}\right)  \partial_{\alpha}V_{N}\left(  Y_{t}^{i,N}-Y_{t}^{j,N}\right)
\left(  Y_{t}^{i,N}-Y_{t}^{j,N}\right)  _{\beta}.
\end{align*}
Let us introduce the function%
\[
\psi_{\alpha\beta}\left(  x\right)  =x_{\beta}\partial_{\alpha}V\left(
x\right)  .
\]
Recall $V_{N}\left(  x\right)  =N^{d}V\left(  Nx\right)  $, hence $\nabla
V_{N}\left(  x\right)  =NN^{d}\left(  \nabla V\right)  \left(  Nx\right)  $,
hence%
\begin{align*}
\partial_{\alpha}V_{N}\left(  Y_{t}^{i,N}-Y_{t}^{j,N}\right)  \left(
Y_{t}^{i,N}-Y_{t}^{j,N}\right)  _{\beta}  &  =NN^{d}\left(  \partial_{\alpha
}V\right)  \left(  N\left(  Y_{t}^{i,N}-Y_{t}^{j,N}\right)  \right)  \left(
Y_{t}^{i,N}-Y_{t}^{j,N}\right)  _{\beta}\\
&  =N^{d}\psi_{\alpha\beta}\left(  N\left(  Y_{t}^{i,N}-Y_{t}^{j,N}\right)
\right)
\end{align*}
therefore%
\[
\left\langle \mu_{t}^{N},\nabla\phi\cdot\nabla\left(  V_{N}\ast\mu_{t}%
^{N}\right)  \right\rangle =\frac{1}{2N^{2d}}\sum_{i,j=1}^{\left[  N^{d}%
\rho\right]  }\sum_{\alpha,\beta=1}^{d}\partial_{\alpha}\partial_{\beta}%
\phi\left(  Y_{t}^{j,N}\right)  N^{d}\psi_{\alpha\beta}\left(  N\left(
Y_{t}^{i,N}-Y_{t}^{j,N}\right)  \right)
\]%
\[
=\frac{1}{2N^{d}}\sum_{i=1}^{\left[  N^{d}\rho\right]  }\sum_{\alpha,\beta
=1}^{d}\partial_{\alpha}\partial_{\beta}\phi\left(  Y_{t}^{j,N}\right)
\sum_{j:\left\vert Y_{t}^{i,N}-Y_{t}^{j,N}\right\vert \leq\frac{r_{1}}{N}}%
\psi_{\alpha\beta}\left(  N\left(  Y_{t}^{i,N}-Y_{t}^{j,N}\right)  \right)  .
\]

It remains to understand where this expression converges.

\subsection{Invariant measures\label{sect inv measure}}

Let us discuss invariant measures for the original microscopic system
(\ref{model in large box}) on the large torus $\mathbb{T}_{N}^{d}$. Let us
parametrize this system in a more general way. Given $\rho>0$ and $L>0$, on
the torus $\mathbb{T}_{L}^{d}=\mathbb{R}^{d}/L\mathbb{Z}^{d}$ (informally
$\left[  0,L\right]  ^{d}$) consider $K_{L}:=\left\lfloor L^{d}\rho
\right\rfloor $ particles $X_{t}^{i}$. The SDE system has the invariant
measure%
\[
\mu_{_{\rho,L}}\left(  dx^{1}...dx^{_{K_{L}}}\right)  =\frac{1}{Z_{_{\rho,L}}%
}\exp\left(  -\frac{1}{2}\sum_{i,j=1}^{K_{L}}V\left(  x^{i}-x^{j}\right)
\right)  dx^{1}...dx^{K_{L}}.
\]

Taken $\psi\in C_{c}^{0}\left(  \mathbb{R}^{d}\right)  $ (it will be one of
the functions $\psi_{\alpha\beta}\left(  x\right)  =x_{\beta}\partial_{\alpha
}V\left(  x\right)  $ introduced above), assume we are able to prove that
there exists $\Psi_{V}:\mathbb{R\rightarrow R}$ such that
\begin{equation}
\lim_{L\rightarrow\infty}\int_{\left(  \mathbb{T}_{L}^{d}\right)  ^{K_{L}}%
}\left\vert \frac{1}{L^{d}}\sum_{i,j=1}^{K_{L}}\psi\left(  x^{i}-x^{j}\right)
-\Psi_{V}\left(  \rho\right)  \right\vert ^{2}\mu_{_{L}}\left(  dx^{1}%
...dx^{_{K_{L}}}\right)  =0. \label{ergodico}%
\end{equation}
If this happens, we could say that ``\textit{spatial averages of local
observables converge}''. This is a very technical result but heuristically
$\frac{1}{L^{d}}\sum_{i,j=1}^{K_{L}}\psi\left(  x^{i}-x^{j}\right)  $ is an
average:
\begin{align*}
\frac{1}{L^{d}}\sum_{i,j=1}^{K_{L}}\psi\left(  x^{i}-x^{j}\right)   &
=\frac{1}{L^{d}}\sum_{i=1}^{K_{L}}F\left(  x\right) \\
F\left(x\right)   &  :=\sum_{\left\vert x^{i}-x^{j}\right\vert \leq
r_{1}}\psi\left(  x^{i}-x^{j}\right)
\end{align*}
with the notation\text{ }$x=\left(  x^{1},...,x^{_{K_{L}}}\right)  $; property
(\ref{ergodico}) is a sort of Law of Large Numbers (or more precisely a
spatial ergodic property) with respect to the Gibbs measure in infinite volume.

\begin{remark}
In discrete systems like those considered in \cite{kipnis1998scaling}, the invariant
measures are much easier, usually product measures and the corresponding
quantities can be computed more explicitly. The weak point of the continuum
theory is the difficulty to compute $\Psi_{V}\left(  \rho\right)  $.
\end{remark}

\subsection{Local equilibrium}

\begin{definition}
\label{def eq loc}Let a continuous function $\rho:\left[  0,T\right]
\times\mathbb{T}^{d}\rightarrow\mathbb{R}$ with constant in time mass
$\overline{\rho}=\int_{\mathbb{T}^{d}}\rho_{t}\left(  y\right)  dy$ be given.
For particle system (\ref{mean field like}) with $K_{N}=\left[  \overline
{\rho}N^{d}\right]  $, we say that local equilibrium holds if, for every
$\varphi,\psi\in C\left(  \mathbb{T}^{d}\right)  $, every $t\in\left[
0,T\right]  $, in probability we have%
\[
\lim_{N\rightarrow\infty}\frac{1}{N^{d}}\sum_{i=1}^{\left[  \overline{\rho
}N^{d}\right]  }\varphi\left(  Y_{t}^{j,N}\right)  \sum_{j:\left\vert
Y_{t}^{i,N}-Y_{t}^{j,N}\right\vert \leq\frac{r_{1}}{N}}\psi\left(  N\left(
Y_{t}^{i,N}-Y_{t}^{j,N}\right)  \right)  =\int_{\mathbb{T}^{d}}\varphi\left(
y\right)  \Psi_{V}\left(  \rho_{t}\left(  y\right)  \right)  dy.
\]

\end{definition}

Assume this property is satisfied. Let us apply it to $\varphi=\partial
_{\alpha}\partial_{\beta}\phi$, $\psi=\psi_{\alpha\beta}$. We get%
\begin{align*}
&  \lim_{N\rightarrow\infty}\frac{1}{2N^{d}}\sum_{i=1}^{\left[  \overline
{\rho}N^{d}\right]  }\partial_{\alpha}\partial_{\beta}\phi\left(  Y_{t}%
^{j,N}\right)  \sum_{j:\left\vert Y_{t}^{i,N}-Y_{t}^{j,N}\right\vert \leq
\frac{r_{1}}{N}}\psi_{\alpha\beta}\left(  N\left(  Y_{t}^{i,N}-Y_{t}%
^{j,N}\right)  \right) \\
&  =\frac{1}{2}\int_{\mathbb{T}^{d}}\partial_{\alpha}\partial_{\beta}%
\phi\left(  x\right)  \Psi_{V}^{\alpha,\beta}\left(  \rho_{t}\left(  x\right)
\right)  dx
\end{align*}
where $\Psi_{V}^{\alpha,\beta}$ is defined from $\psi_{\alpha\beta}$ via
(\ref{ergodico}). hence%
\[
\frac{d}{dt}\left\langle \rho_{t},\phi\right\rangle +\frac{1}{2}%
\int_{\mathbb{T}^{d}}\sum_{\alpha,\beta}\partial_{\alpha}\partial_{\beta}%
\phi\left(  x\right)  \Psi_{V}^{\alpha,\beta}\left(  \rho_{t}\left(  x\right)
\right)  dx=\frac{\sigma^{2}}{2}\left\langle \rho_{t},\Delta\phi\right\rangle
.
\]
By an isotropy argument that we omit we finally get%

\[
\frac{d}{dt}\left\langle \rho_{t},\phi\right\rangle +\frac{1}{2}%
\int_{\mathbb{T}^{d}}\Delta\phi\left(  x\right)  \Psi_{V}\left(  \rho
_{t}\left(  x\right)  \right)  dx=\frac{\sigma^{2}}{2}\left\langle \rho
_{t},\Delta\phi\right\rangle
\]
and by formal integration by parts%
\[
\partial_{t}\rho_{t}=\frac{\sigma^{2}}{2}\Delta\rho_{t}-\frac{1}{2}\Delta
\Psi_{V}\left(  \rho_{t}\right)
\]
where $\Psi_{V}\left(  \rho\right)  $ is $\Psi_{V}^{\alpha,\alpha}\left(
\rho\right)  $ independently of $\alpha=1,...,d$.

\subsection{On local equilibrium:\ a multiscale argument}

It remains to understand the meaning of \ref{def eq loc} and how to prove its
validity. This is well understood in the discrete setting but not so much in
the continuous SDE one. Let us discuss the property
\[
\lim_{N\rightarrow\infty}\frac{1}{N^{d}}\sum_{i,j}\varphi\left(  Y_{t}%
^{j,N}\right)  \psi\left(  N\left(  Y_{t}^{i,N}-Y_{t}^{j,N}\right)  \right)
=\int_{\mathbb{T}^{d}}\varphi\left(  x\right)  \Psi_{V}\left(  \rho_{t}\left(
x\right)  \right)  dx.
\]
Since it will appear several times below, let us write
\[
A_{N,t}:=\frac{1}{N^{d}}\sum_{i,j}\varphi\left(  Y_{t}^{j,N}\right)
\psi\left(  N\left(  Y_{t}^{i,N}-Y_{t}^{j,N}\right)  \right)  .
\]
Decompose $\mathbb{T}^{d}=\left[  0,1\right]  ^{d}$ in $m_{N}$ squares
$Q_{N}\left(  y_{1}^{N}\right)  ,...,Q_{N}\left(  y_{m_{N}}^{N}\right)  $ of
the form%
\[
Q_{N}\left(  y_{k}^{N}\right)  =%
{\displaystyle\prod\limits_{\alpha=1}^{d}}
\left[  y_{k,\alpha}^{N}-\frac{1}{2m_{N}},y_{k,\alpha}^{N}+\frac{1}{2m_{N}%
}\right]  \qquad k=1,...,m_{N}%
\]
where $y_{k}^{N}$ are points of coordinates $y_{k,\alpha}^{N}$, $\alpha
=1,...,d$. The idea is that each one contains a huge number of particles, but
$m_{N}$ is very large as well.\ Hence assume
\begin{align*}
m_{N}  &  \rightarrow\infty\qquad\text{as }N\rightarrow\infty\\
m_{N}  &  =o\left(  N^{d}\right)  .
\end{align*}

Let us suppose the interactions between particles lying in different domains
have a small contribution. Thus let us approximate%
\[
A_{N,t}\sim\frac{1}{m_{N}}\sum_{k=1}^{m_{N}}\frac{m_{N}}{N^{d}}\sum
_{Y_{t}^{i,N},Y_{t}^{j,N}\in Q_{N}\left(  y_{k}^{N}\right)  }\varphi\left(
Y_{t}^{j,N}\right)  \psi\left(  N\left(  Y_{t}^{i,N}-Y_{t}^{j,N}\right)
\right)  .
\]
Since the size $\frac{1}{m_{N}}$ of the little squares is very small,
approximatively%
\[
\varphi\left(  Y_{t}^{j,N}\right)  \sim\varphi\left(  y_{k}^{N}\right)
\qquad\text{when }Y_{t}^{j,N}\in Q_{N}\left(  y_{k}^{N}\right)
\]
and thus let us approximate further%
\[
A_{N,t}\sim\frac{1}{m_{N}}\sum_{k=1}^{m_{N}}\varphi\left(  y_{k}^{N}\right)
\left[  \frac{m_{N}}{N^{d}}\sum_{Y_{t}^{i,N},Y_{t}^{j,N}\in Q_{N}\left(
y_{k}^{N}\right)  }\psi\left(  N\left(  Y_{t}^{i,N}-Y_{t}^{j,N}\right)
\right)  \right]  .
\]
The external average $\frac{1}{m_{N}}\sum_{k=1}^{m_{N}}$ can be seen as the
Riemann sums approximating an integral, yielding%
\[
A_{N,t}\sim\int_{\mathbb{T}^{d}}\varphi\left(  y\right)  \left[  \frac{m_{N}%
}{N^{d}}\sum_{Y_{t}^{i,N},Y_{t}^{j,N}\in Q_{N}\left(  y\right)  }\psi\left(
N\left(  Y_{t}^{i,N}-Y_{t}^{j,N}\right)  \right)  \right]  dy
\]
Thus we have to understand the expression%
\[
\frac{m_{N}}{N^{d}}\sum_{Y_{t}^{i,N},Y_{t}^{j,N}\in Q_{N}\left(  y\right)
}\psi\left(  N\left(  Y_{t}^{i,N}-Y_{t}^{j,N}\right)  \right)
\]
when $N\rightarrow\infty$. Recalling the definition $Y_{t}^{i,N}:=\frac{1}%
{N}X_{N^{2}t}^{i,N}$, the last expression is equal to%
\[
\frac{m_{N}}{N^{d}}\sum_{X_{N^{2}t}^{i,N},X_{N^{2}t}^{j,N}\in%
{\displaystyle\prod\limits_{\alpha=1}^{d}}
\left[  Ny_{\alpha}-\frac{N}{2m_{N}},Ny_{\alpha}+\frac{N}{2m_{N}}\right]
}\psi\left(  X_{N^{2}t}^{i,N}-X_{N^{2}t}^{j,N}\right)  .
\]
Time is much accelerated and thus it is reasonable to approximate $\left\{
X_{N^{2}t}^{i}\right\}  $ with the corresponding stationary process; by the
ergodic property (\ref{ergodico}), the previous expression, computed along
realizations of the stationary process, is close to $\Psi_{V}\left(  \rho
_{t}\left(  x\right)  \right)  $. The choice of the density $\rho_{t}\left(
x\right)  $ is due to the fact that the number of particles in $Q_{N}\left(
y\right)  $ is approximately given by $\rho_{t}\left(  x\right)  \cdot
N/m_{N}$.

\subsection{Summary of the main result}

We have heuristically shown above (proofs are given in  \cite{uchiyama2000pressure},\cite{varadhan1991scaling}) that the
macroscopic limit of the particle system with local interaction is the
nonlinear PDE
\[
\partial_{t}\rho_{t}=\frac{1}{2}\Delta P_{V}\left(  \rho_{t}\right)
\]
where%
\[
P_{V}(\rho):=\sigma^{2}\rho-\Psi_{V}\left(  \rho\right)
\]
and $\Psi_{V}\left(  \rho\right)  $ is given by
\begin{equation}
\Psi_{V}\left(  \rho\right)  =\lim_{L\rightarrow\infty}\frac{1}{L^{d}}%
\sum_{i,j=1}^{K_{L}}\psi_{\alpha\alpha}\left(  x^{i}-x^{j}\right)
\label{ergodico2}%
\end{equation}
for any $\alpha=1,...,d$,
\[
\psi_{\alpha\alpha}\left(  x\right)  =x_{\alpha}\partial_{\alpha}V\left(
x\right)  .
\]
In (\ref{ergodico2}) the sample $x^{1},...,x^{K_{L}}$ must be taken
distributed according to $\mu_{\rho,L}$, see Section \ref{sect inv measure}
and in particular see property (\ref{ergodico}) which explains the meaning of
the limit in (\ref{ergodico2}). In numerical simulations we shall simulate
system (\ref{model in large box}) in a large periodic box and take the tail of
the simulation as a sample approximately distributed as $\mu_{\rho,L}$,
heuristically relying on a result of ergodicity in time.

\section{Quantitative results and conjectures about $P_{V}(\rho)$}\label{sec:quant}

Our aim is to get some information of the function $P_V(\rho)$.  Our idea to estimate $P_V(\rho)$ is to approximate the rhs of \eqref{ergodico} with one of the realization of $(y^1,\cdots,y^N)$ according to the invariant measure $\mu_{\rho,L}$:
\begin{equation}\label{ergodico2_approx}
P_V(\rho)--\sigma^2\rho\approx\lim_{L\rightarrow\infty}\frac{1}{L^{d}}\sum_{i,j=1}^{N_{\rho,L}}\psi\left(
x^{i}-x^{j}\right). 
\end{equation}
with $(x^1,\cdots,x^N)$ realization of the variables $(y^1,\dots,y^N)$, distributed according to $\mu_{\rho,L}$.  At this point it is clear that the key to the study of the virial formula consists in the realization of variables $(y^1,\dots,y^N)$. In the paper we propose two way of visualizing the Gibbs measure: some heuristic arguments followed (and confirmed) by numerical simulations. To produce samples of the Gibbs measure is not a simple task. However, assuming a suitable ergodicity property, we can produce realizations of such measure $\mu^{\rho,L}$ by simulations on large times of a SDE whose invariant measure is $\mu^{\rho,L}$. In a general framework, if $X_t$  denotes the solution of the SDE whose invariant measure is $\mu$, then the time average
\[\frac{1}{T}\int_{t_0}^{t_0+T}\psi(X_s)ds\]
by ergodicity property is a good approximation of the spatial integral
\[\int \psi(x)\mu(dx).\]
Now, applying this idea to our case: the SDE whose invariant measure is $\mu^{\rho,L}$ is given by the system of SDEs
\begin{equation}\label{eq:sdenumerics}
dX^{i,N}=-\sum_{j=1,j\neq i}^{N^{\rho,L}}\nabla U(X^{j,N}_t-X^{i,N}_t)dt+\sigma dB^i_t
\end{equation}
in the torus $\mathbb{T}^d_L$. We perform simulations both in dimensions one and two, even if we focus more on the one dimensional case for a matter of simplicity. 

Both for the heuristic and the numerics we follow the same strategy: we start by making an educated guess about the invariant measure $\mu_{\rho,L}$ depending on the potential $U$ and the density $\rho$. In particular we look for an equilibrium configuration $(\overline{x}^{i,N})_{i=1,\dots,N}$ for a system of deterministic ODEs satisfying 
\begin{equation}\label{eq:equilibriumpoints}
\dot{x}_{t}^{i,N}=-\sum_{j=1,j\neq i}^{N^{\rho,L}}\nabla V(x^{j,N}_t-x^{i,N}_t)
\end{equation}
which the same as system \eqref{eq:sdenumerics} taking $\sigma = 0$. This equilibrium configuration is computed directly for some particular choices of the potential. Of course in some situations multiple equilibrium configurations may exists. In those cases we have to rely on a more intuitive reasoning on what the stationary distribution could be. 
These deterministic equilibrium points are then used for two different purposes. For the numerics, they are used as initial condition for the random dynamics \eqref{eq:sdenumerics}, since we have the intuition that they represents already a good approximation of the invariant measure. Moreover we use the same points as fixed particle positions in the virial formula \eqref{ergodico2_approx} to obtain an explicit expression for $P_{V}(\rho)$. By approximating 
\[
\mu_{\rho,L} \approx \frac{1}{N}\sum_{i=1}^{N} \delta_{\overline{x}^{i,N}}
\]
we can carry out the computational analytically and compare the formula obtained with the values coming from numerical simulations.

In the next sections we presents the results obtained by the two approaches. In particular we analyze the following two main examples for the potential. We will call the first one the \textbf{purely repulsive case}. It is defined by 
\begin{equation}\label{eq:rep_strong}
V(x) = \left(\frac{1}{\alpha|x|^{\alpha}}-C\right) \mathds{1}_{\abs{x}\leq R_{1}},
\end{equation}
where the constant $C$ is chosen in order to make the function $V$ continuous in $|x|=R_1$ and for $\alpha > 0$. It represents the ideal repulsive potential for two solid object that interact by contact forces. Of course what we have in mind is not the case of hard-core interaction between hard spheres, which when touch cannot compenetrate. The example we think of is more that of biological cells, that when touch can stretch an compress a bit one into the other. 
Moreover we also consider an \textbf{attractive-repulsive case}, defined by
\begin{equation}\label{eq:attractiverepulsivepotential}
V(x) = \left(\frac{R_{0}^{\alpha}}{\alpha \abs{x}^{\alpha}}  - \frac{R_{0}^{\beta}}{\beta \abs{x}^{\beta}} +C \right)\mathds{1}_{\abs{x}\leq R_{1}},
\end{equation}
where the constant $C$ is chosen in order to make the function $V$ continuous in $|x|=R_1$, $\alpha > \beta$ and $R_{1} > R_{0}$. This potential has the same shape of the classical Lennard-Jones potential but restricted to have finite range. In this case we target specifically the phenomenon of cell-cell adhesion, where cells repulse each other when too close (due to contact forces)  and tends to stick one to each other when they touch. At longer than contact range there is no interaction, motivating the truncation in the potential $V$.

In both case we  carry out all the analytical computations with more details in dimension one. In higher dimension we only analyze the asymptotic of $P_{V}(\rho)$ and perform numerical simulation to dimension up to two.

\subsection{Purely repulsive case}\label{subsec:repulsive}
\begin{figure}[t]
\centering
\includegraphics[width=0.7\textwidth]{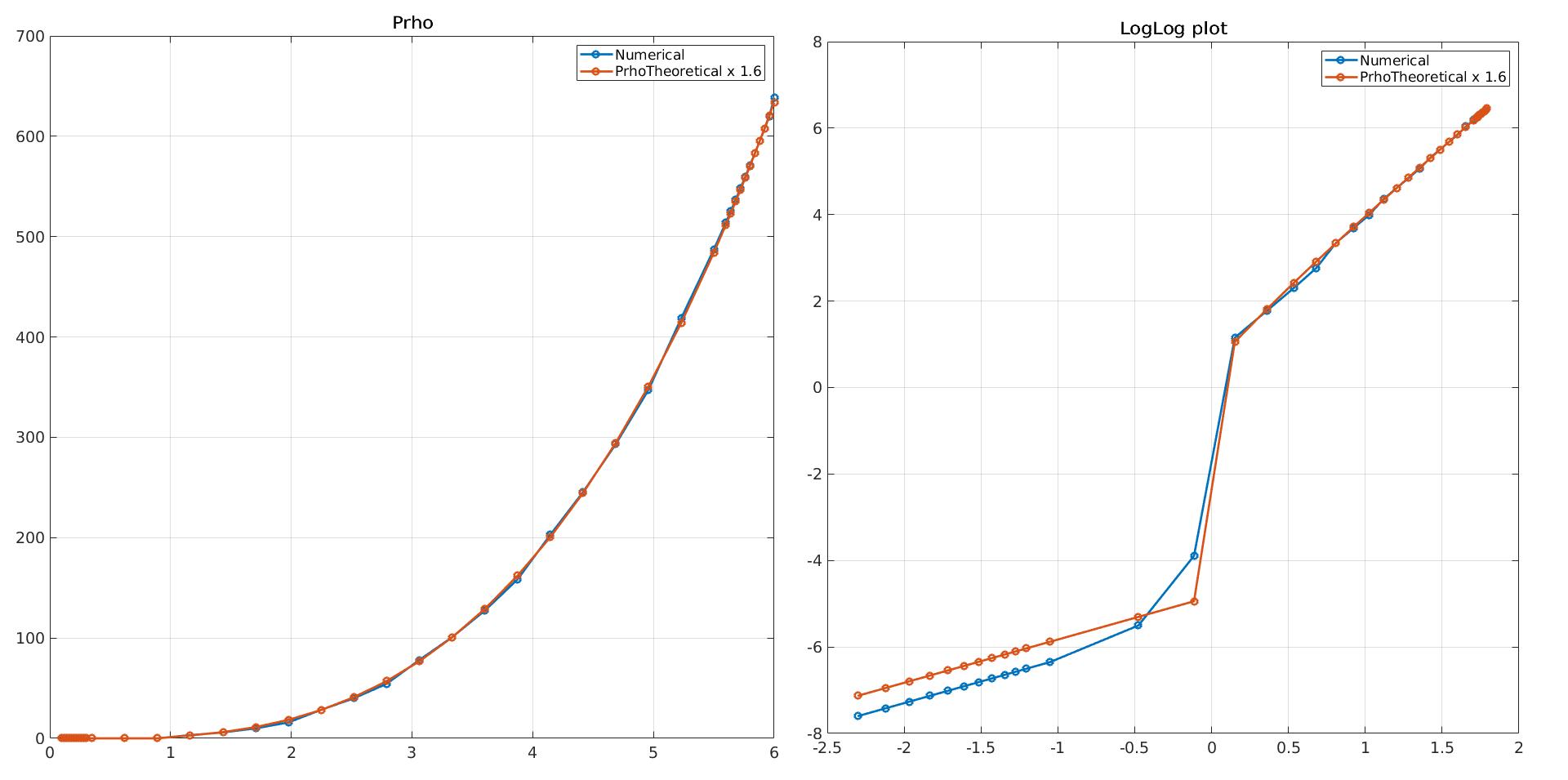}
\caption{Purely repulsive case with compact support, $\alpha = 2,R_{1}=1$ in dimension $d = 1$. Comparison between the function $P_{V}(\rho) $ obtained by numerical simulations and by that obtained in Claim \ref{claim:rep}. Left: comparison between functions $P_{V}(\rho)$ in natural scale. Right: Comparison in $\log\log$ scale. The theoretical function computed as in Claim \ref{claim:rep} has been rescaled by $1.6$.}
\label{fig:rep}
\end{figure}
\begin{figure}[t]
\centering
\vspace{1cm}

\begin{minipage}[b]{0.4\textwidth}
\includegraphics[width=\textwidth]{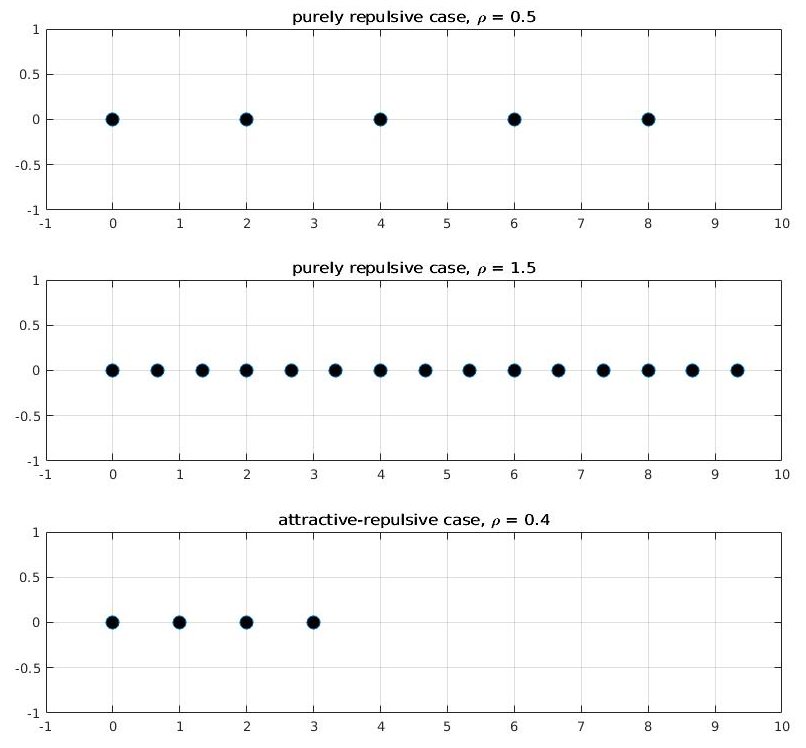}
\end{minipage}\hspace{1.5cm}
\begin{minipage}[t]{0.3\textwidth}
\vspace{-6.5cm}
\begin{tikzpicture}[domain=0.43:2,scale=0.5]
    \draw[->] (-0.5,0) -- (5,0) ; 
    \draw[->] (0,-0.5) -- (0,6);
    \draw[color=red,thick] plot (\x,-0.25+\x^-2);
    \draw[-,color=red,thick] (2,0)--(5,0);    
    \node at (2,-0.5) {$R_1$}; 
\end{tikzpicture} 

\begin{tikzpicture}[domain=0.65:2,scale=0.5]
    \draw[->] (-0.5,0) -- (5,0) ; 
    \draw[->] (0,-0.5) -- (0,6);
    \draw[color=red,thick] plot (\x,+0.234375+\x^-6-2*\x^-3); 
    \draw[-,color=red,thick] (2,0)--(5,0);    
    \node at (1.25,0.3) {\scriptsize $R_0$};  
    \draw[-] (1,0)--(1,-1+0.234375);
\node at (2.2,0.3) {\scriptsize $R_1$};  
\end{tikzpicture} 
\end{minipage}
\caption{Left: Equilibrium configurations in dimension $d = 1$. The distances between consecutive particles is of order $1/\rho$ in the repulsive case, and of order $R_{1}$ in the attractive case for $\rho < 1/R_{1}$. Right: the potential $V(x)$ in the purely repulsive case (top) and attractive repulsive case (bottom).}\label{fig:plottini}
\end{figure}
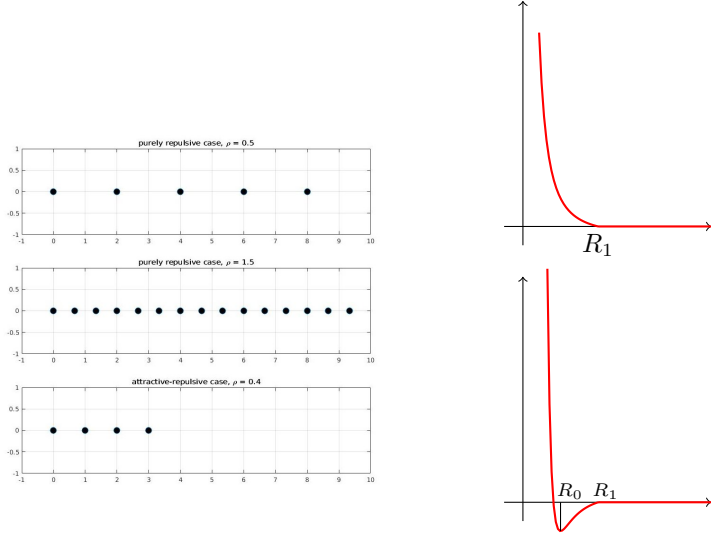

Let us start by the purely repulsive case. We now discuss what are the possible equilibrium points of \eqref{eq:equilibriumpoints} for different values of $\rho$. Let assume that we work in dimension $d=1$ and for simplicity imagine first the case where $R_{1} = +\infty$. Now it is clear that particles have the tendency to stay as far as possible from each other. Assuming a density $\rho$ then the average distance, in a torus of size $L$ is of order $1/\rho$. Hence we approximate 
\[
\mu_{\rho,L} \approx \frac{1}{N}\sum_{i=1}^{N} \delta_{\overline{x}^{i,N}}
\]
with 
\begin{equation}\label{eq:rep_stationary}
   \overline{x}^{i,N} = \frac{i}{\rho}\quad i=1,\dots, N,
\end{equation}
see Figure \ref{fig:plottini}.
If $R_{1} < \infty$, hence the potential has compact support, in order for particle to interact with each other one needs to have $1/\rho < R_{1}$, in other terms $\rho > 1/R_{1}$. We figure out that the behavior of particles is not so different from the long range case. Let us focus on the low density configuration $\rho<1/R_1$ and let us assume particles to be at equal distance $1/\rho$ . We can think that each Brownian particle $X^i$ has attached to itself a sort of delimitation zone of amplitude $R_1$, $[X^i-R_1/2,X^i+R_1/2]$ and as soon another particle goes inside such area, it is pushed away. So at the end particles will always come back to the equal distance configuration. For $\rho>1/R_1$ we expect particle to stay as large as possible. So also for the compactly supported case we assume that \eqref{eq:rep_stationary} holds.

\begin{claim}
\label{claim:rep}
In dimension $d = 1$, if the potential $U$ is given by \eqref{eq:rep_strong} we claim that 
\[P_V(\rho)\approx\begin{cases}
\sigma^2\rho &\rho<1/R_1,\\
\sigma^{2}\rho  +\frac{\alpha}{\alpha-1}\rho^{1+\alpha}- \frac{1}{R_{1}^{\alpha-1}(\alpha-1)}\rho^{2} &\rho\geq 1/R_1
\end{cases}\]
for $\alpha \neq 1$.
Instead for $\alpha = 1$ we claim
\[P_V(\rho)\approx\begin{cases}
\sigma^2\rho &\rho<1/R_1,\\
\sigma^{2}\rho  +\rho^{2} (\log(\rho) + \log(R_{1})+1) &\rho\geq 1/R_1.
\end{cases}\]
\end{claim}
Let us start by discussing the case $\rho < 1/R_{1}$. By using equation \eqref{ergodico2_approx} with points $\overline{x}^{i,N}$ we obtain
\begin{equation*}
-\frac{1}{L}\sum_{i=1}^{N_{\rho,L}}\sum_{j=1}^{N_{\rho,L}}V'\left(\abs{\overline{x}^{i,N}-\overline{x}^{j,N}}\right)\abs{\overline{x}^{i,N}-\overline{x}^{j,N}} = -\rho \sum_{i=1}^{\rho L}V'\left(\abs{\overline{x}^{i,N}}\right)\abs{\overline{x}^{i,N}}\\
\end{equation*}
since $N_{\rho L} = \rho L$ and particles are all at the same distance one from the other. Moreover we have
\[
-\rho \sum_{i=1}^{\rho L}V'\left(\abs{\overline{x}^{i,N}}\right)\abs{\overline{x}^{i,N}} = \rho\sum_{i=1}^{\rho L} \frac{\rho^{\alpha+1}}{i^{\alpha+1}} \frac{i}{\rho}\mathds{1}_{\abs{\frac{i}{\rho} } \leq R_{1}} = 0
\]
since we assumed $1/\rho > R_{1}$. In the case $\rho > 1/R_{1}$ we can carry out the computation (assume $\alpha \neq 1$ for simplicity) obtaining
\[
\rho\sum_{i=1}^{\rho L} \frac{\rho^{\alpha+1}}{i^{\alpha+1}} \frac{i}{\rho}\mathds{1}_{\abs{\frac{i}{\rho} } \leq R_{1}} = \rho^{1+\alpha}\sum_{i=1	}^{\rho R_{1}} \frac{1}{i^{\alpha}}.
\]
Now we approximate
\[
\sum_{i=1	}^{\rho R_{1}} \frac{1}{i^{\alpha}} \approx 1+\int_{1}^{\rho R_{1}}\frac{1}{x^{\alpha}} \,dx = 1+\frac{1}{\alpha-1} \left( 1-	\frac{1}{R_{1}^{\alpha-1}\rho^{\alpha-1}} 	 \right),
\]
hence
\[
\rho^{1+\alpha}\sum_{i=1	}^{\rho R_{1}} \frac{1}{i^{\alpha}} = \rho^{1+\alpha} \left[1+\frac{1}{\alpha-1} \left( 1-	\frac{1}{R_{1}^{\alpha-1}\rho^{\alpha-1}} 	 \right)  \right]\\
= \frac{\alpha}{\alpha-1}\rho^{1+\alpha}- \frac{1}{R_{1}^{\alpha-1}(\alpha-1)}\rho^{2}.
\]
The case where $\alpha = 1$ is derived in the same manner, by a proper estimation of the sum $\sum_{i=1	}^{\rho R_{1}} \frac{1}{i}$.

Let us point out that the asymptotic for large values of $\rho$ coincide with what is expressed in Lemma 8.5 in \cite{uchiyama2000pressure}. In fact for $\alpha > 1$ i.e. if $U(r)$ is not in $L^{1}$ around zero, the leading term in the expression presented is of order $\rho^{1+\alpha}$. Moreover, we also observe that for those values of $\alpha$ that makes $U(r)$ integrable around zero ($\alpha < 1$) the leading order of growth is again $\rho^{2}$. This result is in agreement with what has been proven in \cite{uchiyama2000pressure} for integrable potentials and it also confirms the order of growth derived by Oelschläger in \cite{oelschlager1990large}. \\
In summary our claim is that: for large values of $\rho$ the asymptotic of $P_{V}(\rho)$ obeys to the following rule:
\begin{equation}
\begin{cases}
\rho^{1+\alpha} &\text{ if }\alpha > 1,\\
\rho^{2} &\text{ if } \alpha < 1,\\
\rho^{2}\log(\rho) &\text{ if } \alpha = 1.
\end{cases}
\end{equation}
For small values of $\rho$ the behavior instead is fully linear in all the analyzed cases.\\
The asymptotic for small and large $\rho$ is confirmed by numerical simulations. In fact we computed the function $P_{V}(\rho)$ for different values of $\rho$, very small for the small density asymptotic and very large for the high density. Afterwards we computed an approximation of the growth exponent by measuring the slope of the function $P_{V}(\rho)$ in $\log\log$-scale. All the results are presented in Table \ref{tab:repulsive}.
\begin{table}[t]
\centering
\scalebox{0.70}{
\begin{tabular}{rcccc}
\multicolumn{1}{c}{}                & \multicolumn{2}{c}{$\rho \approx 0.2$}                                                                                       & \multicolumn{2}{c}{$\rho \approx 6$}                                                                    \\
\multicolumn{1}{l|}{}               & $V(x) = \left(\frac{1}{\alpha\abs{x}^{\alpha}}-C\right)$ & \multicolumn{1}{c|}{$V(x) = \left(\frac{1}{\alpha\abs{x}^{\alpha}}-C\right)\mathds{1}_{\abs{x}\leq R_{1}}$} & $V(x) = \left(\frac{1}{\alpha\abs{x}^{\alpha}}-C\right)$ & $V(x) = \left(\frac{1}{\alpha\abs{x}^{\alpha}}-C\right)\mathds{1}_{\abs{x}\leq R_{1}}$ \\ \hline
\multicolumn{1}{r|}{$\alpha = 0.5$} & $1.001$                             & \multicolumn{1}{c|}{$1.013$}                                                           & $2.04$                              & $1.48$                                                            \\
\multicolumn{1}{r|}{$\alpha = 1$}   & $1.03$                              & \multicolumn{1}{c|}{$1.025$}                                                           & $2.13$                              & $2.03$                                                            \\
\multicolumn{1}{r|}{$\alpha = 2$}   & $1.027$                             & \multicolumn{1}{c|}{$1.017$}                                                           & $3.004$                             & $3.004$                                                           \\
\multicolumn{1}{r|}{$\alpha = 3$}   & $1.025$                             & \multicolumn{1}{c|}{$1.023$}                                                           & $4.007$                             & $4.007$                                                           \\
\multicolumn{1}{r|}{$\alpha = 4$}   & $1.009$                             & \multicolumn{1}{c|}{$1.011$}                                                           & $5.003$                             & $5.003$                                                          
\end{tabular}}
\caption{ Leading exponent for the function $P_V(\rho)$ for various potential, for different values of $\alpha$ and $\rho$ in dimension $d = 1$. Note that there are some discrepancies in the case $\alpha < 1$ and $\rho \approx 6$ since one recover the theoretical behavior proved in [13] only in the case where the potential is not compactly supported. }
\label{tab:repulsive}
\end{table}

Moreover we also compare directly the values obtained by the function $P_{V}(\rho)$ computed by numerical simulations and that obtained by the formula of Claim \ref{claim:rep}. This comparison is highlighted in Figure \ref{fig:rep}. As we can see the global regime is captured nicely by the function captured by heuristics. We also compare the plots in $\log\log$ scale, in order to have a better insight on what is the leading term for each values of $\rho$. In both plots the function proposed in Claim \ref{claim:rep} has been multiplied by a coefficient $1.6$ to obtain a better approximation. This shows that the conjecture proposed in Claim \ref{claim:rep}, even if correct in terms of the growth exponent, but still needs a small refining in term of multiplicative constants.

\begin{center}
\begin{figure}[]
\begin{minipage}[t]{0.4\textwidth}
\includegraphics[width=\textwidth]{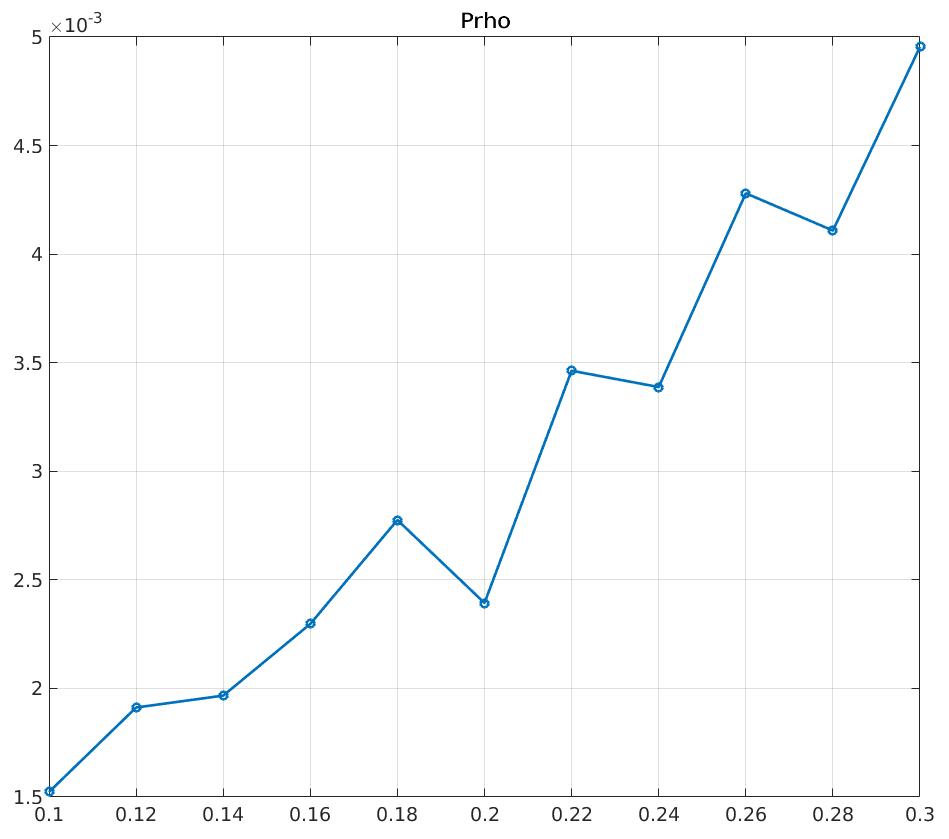}
\end{minipage}\hspace{0.5cm}
\begin{minipage}[t]{0.1\textwidth}\vspace{-4cm}
\scalebox{0.80}{
\begin{tabular}{rc}
\multicolumn{1}{c}{}                           & $\rho \approx 6$                                                                                                                                  \\
\multicolumn{1}{l|}{}                          & $V(x) = \left(\frac{1}{\alpha \abs{x}^{\alpha}}  - \frac{1}{\beta \abs{x}^{\beta}} +C \right)\mathds{1}_{\abs{x}\leq 1.5}$ \\ \hline
\multicolumn{1}{r|}{$\alpha = 2, \beta = 1.5$} & $3.01$                                                                                                                                              \\
\multicolumn{1}{r|}{$\alpha = 3, \beta = 2$}   & $4.03$                                                                                                                                              \\
\multicolumn{1}{r|}{$\alpha = 4, \beta = 3$}   & $5.007$                                                                                                                                            
\end{tabular}}
\end{minipage}\caption{Left: Plots of the function $P_{V}(\rho)$ in the attractive-repulsive case computed numerically, for small values of $\rho$. The estimated exponent by $\log\log$ scale is $1.0265$, confirming the linearity expressed in \eqref{eq:Prho_h}. The slope estimated is around $0.0175$ which does not coincide with the vale expected $\sigma^{2}/2 = 0.005$. Right: Leading exponent of the function $P_{V}(\rho)$ in the attractive-repulsive case for different values of $\alpha,\beta$ and large $\rho$ in dimension $d = 1$.} \label{fig:small+tab}
\end{figure}
\end{center}

In the multidimensional case the situation is more complicate, therefore we only sketch the main computation to derive the leading order term for high values of $\rho$. Assume dimension $d$ and assume $N = \rho L^{d}$ to be such that $N^{1/d}$ is an integer. In the multidimensional case is less clear how to compute equilibrium points $(\overline{x}^{i,N})_{i=1,\dots,N}$. Hence for simplicity we assume that points $\overline{x}^{i,N})$ are arranged on a uniform grid $\mathbb{Z}^{d}/L\mathbb{Z}^{d}$ embedded inside $\mathbb{T}_{L}^{d}$:
\[
\overline{x}^{i,N} = \left(\frac{i_{1}}{\rho^{1/d}},\dots, \frac{i_{d}}{\rho^{1/d}} \right) \quad  i_{k} = 1,\dots,N^{1/d}.
\]
With this choice the distance of the $i$-th particle and the particle positioned in the origin is
\[
\frac{\abs{(i_{1},\dots,i_{d})}}{\rho^{1/d}}.
\]
By formula \eqref{ergodico2_approx} we obtain
\[
-\frac{1}{L^{d}}\sum_{i=1}^{N_{\rho,L}}\sum_{j=1}^{N_{\rho,L}}V'\left(\abs{\overline{x}^{i,N}-\overline{x}^{j,N}}\right)\abs{\overline{x}^{i,N}-\overline{x}^{j,N}} = -\rho \sum_{i=1}^{\rho L^{d}}V'\left(\abs{\overline{x}^{i,N}}\right)\abs{\overline{x}^{i,N}}
\]
\[
= \rho \sum_{i=1}^{\rho L^{d}} \frac{\rho^{\alpha/d}}{\abs{(i_{1},\dots,i_{d})}^{\alpha}} \mathds{1}_{\abs{(i_{1},\dots,i_{d})}\leq R_{1}\rho^{1/d}} = \rho^{1+\alpha/d} \sum_{i=1}^{\rho L^{d}} \frac{1}{\abs{(i_{1},\dots,i_{d})}^{\alpha}} \mathds{1}_{\abs{(i_{1},\dots,i_{d})}\leq R_{1}\rho^{1/d}}.
\]
Carrying out an explicit computation here is more difficult. However in the case where $R_{1} = +\infty$ and $\alpha > d$, since  the series 
\[
\sum_{i=1}^{\infty} \frac{1}{\abs{(i_{1},\dots,i_{d})}^{\alpha}} < \infty
\]
we understand that the leading term in the expression of $P_{V}(\rho)$ is $\rho^{1+\alpha/d}$. This intuition is confirmed by the values computed numerically in the case $d = 2$, Table \ref{tab:repulsive d=2}. 
\begin{table}
\centering
\scalebox{0.9}{
\begin{tabular}{rcc}
\multicolumn{1}{c}{}              & \multicolumn{2}{c}{$\rho \approx 6$}                                                                                       \\
\multicolumn{1}{l|}{}             & $V(x) = \left(\frac{1}{\abs{x}^{\alpha}}-C\right)$ & \multicolumn{1}{c|}{$V(x) =\left( \frac{1}{\abs{x}^{\alpha}}-C\right)\mathds{1}_{\abs{x}\leq R_{1}}$} \\ \hline
\multicolumn{1}{r|}{$\alpha = 1$} & $2.031$                             & \multicolumn{1}{c|}{$1.421$}                                                           \\
\multicolumn{1}{r|}{$\alpha = 2$} & $2.11$                              & \multicolumn{1}{c|}{$2.04$}                                                            \\
\multicolumn{1}{r|}{$\alpha = 3$} & $2.503$                             & \multicolumn{1}{c|}{$2.511$}                                                           \\
\multicolumn{1}{r|}{$\alpha = 4$} & $3.04$                              & \multicolumn{1}{c|}{$3.07$}                                                            \\
\multicolumn{1}{l|}{$\alpha = 5$} & $3.54$                              & \multicolumn{1}{c|}{$3.55$}                                                           
\end{tabular}
}
\caption{Leading exponent for the function $P_{V}(\rho)$ for various potential, for different values of $\alpha$ for large $\rho$ in dimension $d = 2$.}
\label{tab:repulsive d=2}
\end{table}

\subsection{Attractive-repulsive potential}\label{subsec:attractive_repulsive}

We now consider the attractive-repulsive case, see Potential given by \eqref{eq:attractiverepulsivepotential}. In this case the potential $U(x)$ change the sign, see Figure \ref{fig:plottini}. Let us now start to discuss equilibrium configurations for particles \eqref{eq:equilibriumpoints} in this case. Again we propose the argument in dimension $d=1$ for a matter of simplicity.
In this case the minimum of the potential $U$ is obtained in the values $R_{0}$. Therefore particles have the tendency to align at a distance of $R_{0}$ one to each other. Moreover we will always analyze the case where $R_{1} < 2R_{0}$. This conditions specifies the fact that, if a particle is located at distance $R_{0}$ from its two neighbors, then it interacts only with the two of them, since other particles are outside the support of the potential $V(x)$. If we think particles as cells of radius $R_0/2$, we are describing the phenomenon of cellular adhesion, indeed attraction occurs only when cells are in contact, namely when they are at distance smaller then $2R_0$. This assumption creates an uniform structure: the equilibrium configuration consists of particles collocated at equal distance $R_{0}$, see Figure \ref{fig:plottini}.
 It is worth noting that this configuration is not uniformly spread in the whole torus. In fact for low values of $\rho$ particles will tend to form a cluster that occupies only a portion of the available space. 
When the density is higher, $\rho > 1/R_{0}$, particles are compressed to be closer than $R_{0}$ one from the other, and we recover the uniform distribution at distance $1/\rho$ as in the purely repulsive case. This intuition is confirmed by numerical simulations where we had the system start at a random initial condition and let it evolve towards its stationary configuration.
Summarizing we imagine the equilibrium configuration in the attractive-repulsive case to be 
\[
\begin{cases}
\text{equally spaced particles at distance $R_{0}$} &\text{ if } \rho < 1/R_{0},\\
\text{equally spaced particles at distance $1/\rho$} &\text{ if } \rho \geq 1/R_{0}.\\
\end{cases}
\]
Notice that the situation here is different from the purely repulsive case particles will always place at a distance $1/\rho$ for every value of $\rho$.
\begin{figure}[h!]
\centering
\includegraphics[width=0.7\textwidth]{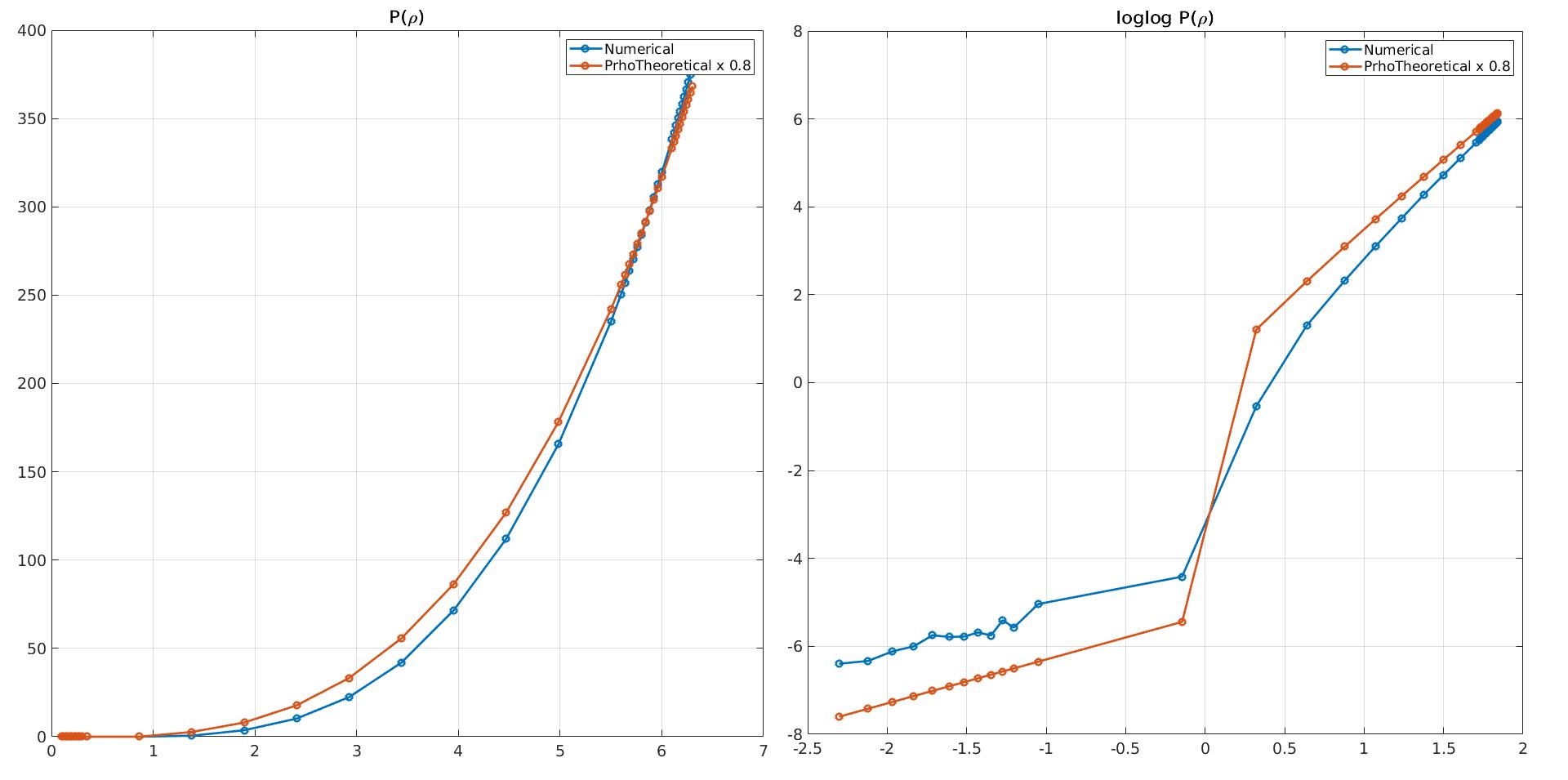}
\caption{Attractive repulsive case with compact support, $\alpha = 2,\beta = 1.5,R_{1}=1.5$ in dimension $d = 1$. Comparison between the function $P_{V}(\rho) $ obtained by numerical simulations and by that obtained by heuristics. Left: comparison between functions $P_{V}(\rho)$ in natural scale. Right: Comparison in $\log\log$ scale. The theoretical function computed as in Claim \ref{claim:rep} has been rescaled by $0.8$.}
\label{fig:attr}
\end{figure}
We now present a claim for the function $P_{V}(\rho)$ in the attractive-repulsive case. For simplicity we will present only the case where $\alpha > \beta > 1$. 
\begin{claim}\label{claim:attrrep}
In dimension $d = 1$, if the potential $U$ is given by \eqref{eq:attractiverepulsivepotential} we claim that 
\[P_V(\rho)\approx\begin{cases}
\sigma^2\rho +O(\rho)&\rho<1/R_0,\\
\sigma^{2}\rho  +\frac{\alpha}{\alpha-1}\rho^{1+\alpha} - \frac{\beta}{\beta-1}\rho^{1+\beta} - \left(\frac{1}{R_{1}^{\alpha-1}(\alpha-1)} - \frac{1}{R_{1}^{\beta-1}(\beta-1)}\right) \rho^{2}&\rho\geq 1/R_{0}
\end{cases}\]
\end{claim}
\noindent By repeating the same argument proposed on the repulsive case, we get as approximation for $P_V(\rho)$:
\begin{equation}\label{eq:Prho_h}
P^{h}_V(\rho)\approx\begin{cases}
\sigma^2\rho&\rho<1/R_0,\\
\sigma^{2}\rho  +\frac{\alpha}{\alpha-1}\rho^{1+\alpha} - \frac{\beta}{\beta-1}\rho^{1+\beta} - \left(\frac{1}{R_{1}^{\alpha-1}(\alpha-1)} - \frac{1}{R_{1}^{\beta-1}(\beta-1)}\right) \rho^{2}&\rho\geq 1/R_{0}
\end{cases}
\end{equation}
In this case numerical simulations and the heuristic result show some tiny discrepancies. We observe some difference regarding the multiplicative constant but no differences are observed on the side of growth exponent. First we analyze the behavior for small $\rho$, Figure \ref{fig:small+tab}. Here we see that, while the growth is indeed linear, the slope of the curve does not coincide with the value $\sigma^{2}/2$ which has been conjectured in \eqref{eq:Prho_h}. This shows that in the attractive repulsive case there could be some additional contribution of the order $O(\rho)$ that has not yet been identified. Moreover we also analyze, via the slope of the function $P_{V}(\rho)$ in $\log\log$ scale, the leading exponent for large values of $\rho$. The results are collected again on Figure \ref{fig:small+tab}. Here we confirm that the leading asymptotic exponent is given by $1+\alpha$ as in the purely repulsive case. In summary we can conclude that up to some correction term, not clearly identified, but of order $O(\rho)$ the function proposed by the heuristic argument catches the global behavior of $P_V(\rho)$.


\end{document}